\documentclass[10pt,reqno]{amsart}
\usepackage{amsmath}
\usepackage{amssymb}
\usepackage{amsthm}
\usepackage{xcolor}

\newlength{\defbaselineskip}
\newcommand{\setlinespacing}[1]%
           {\setlength{\baselineskip}{#1 \defbaselineskip}}

\numberwithin{equation}{section}

\newtheorem{thm}{Theorem}[section]

\theoremstyle{definition}

\theoremstyle{remark}
\newtheorem{rem}[thm]{Remark}
\numberwithin{equation}{section}

\begin{document}

\title[Reverse H\"older inequality]
{On a reverse H\"older inequality for Schr\"odinger operators}

\author{Seongyeon Kim and Ihyeok Seo}

\thanks{This research was supported by a KIAS Individual Grant (MG082901) at Korea Institute for Advanced Study (S. Kim) and NRF-2019R1F1A1061316 (I. Seo).}


\address{School of Mathematics, Korea Institute for Advanced Study, Seoul 02455, Republic of Korea}
\email{synkim@kias.re.kr}

\address{Department of Mathematics, Sungkyunkwan University, Suwon 16419, Republic of Korea}
\email{ihseo@skku.edu}

\begin{abstract}
We obtain a reverse H\"older inequality for the eigenfuctions of the Schr\"odinger operator
with slowly decaying potentials.
The class of potentials includes singular potentials which decay like $|x|^{-\alpha}$ with $0<\alpha<2$, in particular the Coulomb potential.
\end{abstract}

\maketitle

\section{Introduction}
In this paper we are concerned with a reverse H\"older inequality for the eigenfunctions of
the Schr\"odinger operator $-\Delta+V(x)$ in $L^2(\mathbb{R}^n)$.
More generally, we consider second-order elliptic operators of the form
\begin{equation*}
    Lu=- \sum_{i,j=1}^n \frac{\partial}{\partial x_j} (a_{ij}(x)u_{x_i}) + V(x)u
 \end{equation*}
where $a_{ij}(x)$ is a measurable and real-valued function, and the matrix $(a_{ij}(x))_{n\times n}$ is uniformly elliptic.
Namely, there exists a positive constant $\Lambda$ such that
\begin{equation}\label{ellipticity}
    \sum_{i,j=1}^n a_{ij}(x) \xi_i\xi_j \geq \Lambda |\xi|^2
\end{equation}
for $x,\xi \in \mathbb{R}^n $.
Particularly when $a_{ij}=\delta_{ij}$ (Kronecker delta function),
the operator $L$ becomes equivalent to the classical Schr\"odinger operator.
In this regard, we shall call a real-valued function $V(x)$ the potential.

Reverse H\"older inequalities for solutions to the following Dirichlet boundary problem have been studied for a long time:
\begin{equation}\label{maineq}
\begin{cases}
\begin{split}
Lu=- \sum_{i,j=1}^n \frac{\partial}{\partial x_j} (a_{ij}(x)u_{x_i}) + V(x)u &= \lambda u \,\,\quad \text{in} \,\,\,\,\Omega \\
u &=0 \qquad \text{on} \,\,\,\partial \Omega
\end{split}
\end{cases}
\end{equation}
where $\Omega$ is a bounded region in $ \mathbb{R}^{n}$.
When $n=2$, Payne and Rayner \cite{PR} showed that if $\lambda$ is the first eigenvalue and $u$ is the corresponding eigenfunction
of the problem \eqref{maineq} with $a_{ij}=\delta_{ij}$ and $V(x)\equiv0$,
\begin{equation*}
\begin{cases}
\begin{split}
- \Delta u &= \lambda u \,\,\quad \text{in} \,\,\,\,\Omega \\
u &=0 \qquad \text{on} \,\,\,\partial \Omega,
\end{split}
\end{cases}
\end{equation*}
then the following reverse Schwarz inequality holds:
\begin{equation*}
   \|u\|_{L^2(\Omega)} \leq \sqrt{\frac{\lambda}{4 \pi}} \|u\|_{L^1(\Omega)}.
\end{equation*}
This result was extended to higher dimensions by Kohler-Jobin \cite{K} (see also \cite{PR2}).
In the general setting \eqref{maineq}, the reverse H\"older inequalities,
\begin{equation*}
\|u\|_{L^q(\Omega)} \leq C_{p,q,\lambda,n} \|u\|_{L^p(\Omega)},\quad q\geq p>0,
\end{equation*}
were obtained later by
Talenti \cite{T} for $q=2$ and $p=1$, and by Chiti \cite{Ch} for all $q\geq p>0$, but with a nonnegative potential $V\geq0$
and with symmetric coefficients $a_{ij}=a_{ji}$.

Our aim in this paper is to remove these restrictions. Namely, we obtain a reverse H\"older inequality for solutions of \eqref{maineq}
where $V$ is allowed to be negative and we do not need to assume the symmetry, $a_{ij}=a_{ji}$.

For the purpose, we use a completely different approach based on a combination between the Fefferman-Phong inequality and the classical Moser's iteration technique.
Compared with the approach, the authors of the previous results \cite{PR,K,Ch} mentioned above are mainly interested on isoperimetric inequalities by which they obtain explicit constants in their reverse H\"older inequalities and could characterize equality which is not the main issue in the present work.

Before stating our results, we introduce the Morrey-Campanato class ${\mathcal{L}}^{\alpha,r}$ of potentials $V$,  which is defined for $\alpha>0$ and $1\leq r\leq n/\alpha$ by
\begin{equation*}
V \in {\mathcal{L}}^{\alpha,r} \quad \Leftrightarrow \quad \sup_{x \in \mathbb{R}^n, \rho>0}\rho^{\alpha-n/r} \left( \int_{B_\rho(x)} |V(y)|^r dy \right)^{1/r} < \infty,
\end{equation*}
where $B_\rho(x)$ is the ball centered at $x$ with radius $\rho$.
In particular, ${\mathcal{L}}^{\alpha,n/\alpha} = L^{n/\alpha}$ and $1/|x|^\alpha \in L^{n/\alpha, \infty} \subset {\mathcal{L}}^{\alpha,r}$
if $1 \le r < n/\alpha$.
Let us next make precise what we mean by a weak solution of the problem \eqref{maineq}.
We say that a function $u \in H_0^1(\Omega)$ is a weak solution if
\begin{equation}\label{idde}
\int_{\Omega} \sum_{i,j=1}^{n} a_{ij}(x) u_{x_i} \phi_{x_j} dx+ \int_{\Omega}V(x) u(x) \phi (x) dx =  \int_{\Omega}\lambda u(x) \phi(x) dx
\end{equation}
for every $\phi \in H_0^1(\Omega)$.
Our result is then the following theorem.

\begin{thm}\label{thm1}
Let $n\geq3$. Assume that $ u \in H_{0}^1 (\Omega) $ is a weak solution of the problem \eqref{maineq}
with $\lambda\in\mathbb{R}$ and $V\in\mathcal{L}^{\alpha,r}$ for $\alpha<2$ and $r>2/\alpha$.
Then we have
\begin{equation}\label{RHI}
\|u\|_{L^{q} (\Omega)}\leq CC_\alpha^{\frac{n}{2p}} \max\{p,2\}^{\frac{n}{p(2-\alpha)}} \Big(\frac{n}{n-2}\Big)^{\frac{n(n-2)}{p(2-\alpha)}}\|u\|_{L^p (\Omega)}
\end{equation}
for all $q\geq p>0$. Here, $C$ is a constant depending on $\Lambda, \lambda, p, q, n$ and $\Omega$,
and
$$C_\alpha=1+\alpha^{\frac{\alpha}{2-\alpha}}\Big(\frac{2C_n}{\Lambda}\|V\|_{\mathcal{L}^{\alpha, r}}\Big)^{2/(2-\alpha)}$$
with a constant $C_n$ depending on $n$ and arising from the Fefferman-Phong inequality \eqref{FP}.
\end{thm}

\begin{rem}
The class $\mathcal{L}^{\alpha,r}$, $\alpha<2$, of potentials in the theorem includes the positive homogeneous potentials $a|x|^{-\alpha}$ with $a>0$ and $0<\alpha<2$ in three and higher dimensions, in particular the Coulomb potential.
\end{rem}

In Section \ref{sec2} we prove Theorem \ref{thm1}.
Compared with the previous results \cite{K,Ch} based on rearrangements of functions, our approach works also for negative potentials and for non-symmetric coefficients $a_{ij}$.
We use a completely different approach based on a combination between the Fefferman-Phong inequality and the classical Moser's iteration technique.
By making use of a two-dimensional analogue of the Fefferman-Phong inequality, we obtain a reverse H\"older inequality for $n=2$ as well.
See Section \ref{sec3} for details.

Throughout this paper, we denote $A \lesssim B$ to mean $A \leq CB$ with unspecified constant $C>0$ which may be different at each occurrence.

\section{Proof of Theorem \ref{thm1}}\label{sec2}

In this section we prove the reverse H\"older inequality \eqref{RHI}.
Since a complex-valued solution $u$ satisfies \eqref{idde} for every complex $\phi\in H_0^1(\Omega)$,
one can easily see that real and imaginary parts of the solution also satisfy \eqref{idde} for every real $\phi\in H_0^1(\Omega)$.
On the other hand, once we prove the inequality for the real and imaginary parts, we get the same inequality for $u$. Indeed, using the inequality $(a+b)^s \leq C(a^s + b^s)$ for $a,b>0$ and $s>0$, one can see
\begin{align*}
\|u\|_{L^q (\Omega)}^q &= \int_{\Omega} (|\textrm{Re}\,u+i\textrm{Im}\,u|^2)^{q/2} dx  \\
&= \int_{\Omega} (|\textrm{Re}\,u|^2 + |\textrm{Im}\,u|^2 )^{q/2} dx  \\
&\leq C \int_{\Omega} |\textrm{Re}\,u|^q + |\textrm{Im}\,u|^q dx \\
&\leq C(\|\textrm{Re}\,u\|_{L^q (\Omega)}^q + \|\textrm{Im}\,u\|_{L^q (\Omega)}^q) \\
&\leq C(\|\textrm{Re}\,u\|_{L^p (\Omega)}^q + \|\textrm{Im}\,u\|_{L^p (\Omega)}^q) \\
&\leq C\|u\|_{L^p (\Omega)}^q.
\end{align*}
Hence we may assume that the solution $u$ is a real-valued function.

Now we decompose $u$ into two parts, $f= \max \{u,0 \}$ and $g=-\min\{u,0\}$.
Then it is enough to prove that \eqref{RHI} holds for $f$ and $g$.
Indeed,
\begin{align*}
\|u\|_{L^{q}(\Omega)} &= \|f - g\|_{L^{q}(\Omega)} \\
&\leq \|f\|_{L^{q}(\Omega)} + \|g\|_{L^{q}(\Omega)} \\
&\leq C(\|f\|_{L^{p}(\Omega)} + \|g\|_{L^{p}(\Omega)}) \\
&\leq C\|u\|_{L^{p}(\Omega)}.
\end{align*}
We only consider $f$ because the proof for $g$ follows obviously from the same argument.
To prove \eqref{RHI} for $f$, we now divide cases into two parts, $p\geq 2$ and $p<2$.

\subsection{The case $p \ge 2$}
In this case we will show that for all $\tau \geq 2$
\begin{equation}\label{ineq1}
\|f\|_{L^{\tau \omega} (\Omega)} \lesssim C_\alpha^{1/\tau} \tau^{\frac{2}{\tau(2-\alpha)}} \|f\|_{L^{\tau} (\Omega)}
\end{equation}
with $\omega = n/(n-2)$.
Beginning with $\tau = p$, we then iterate as $\tau= p, p\omega , p\omega^2 , \cdots,$ to obtain \eqref{RHI}.
Indeed, first put
\begin{equation*}
\tau_i = p \omega^i
\end{equation*}
for $i=0,1,2,\cdots$.
Since $\tau_i = \tau_{i-1} \omega$, we then get
\begin{align*}
\|f\|_{L^{\tau_i} (\Omega)} &\lesssim C_\alpha^{\frac1{\tau_{i-1}}} {\tau_{i-1}}^{\frac{2}{(2-\alpha)\tau_{i-1}}} \|f\|_{L^{\tau_{i-1}}(\Omega)} \\
&= C_\alpha^{\frac1{p\omega^{i-1}}} \left( p\omega^{i-1} \right)^{\frac2{p(2-\alpha)\omega^{i-1}}} \|f\|_{L^{\tau_{i-1}}(\Omega)}
\end{align*}
for $i=1,2,\cdots,$ which implies by iteration that
\begin{equation*}
\|f\|_{L^{\tau_i}(\Omega)}
\lesssim \big( C_\alpha p^{{2}/(2-\alpha)} \big)^{\sum_{k=1}^i (p\omega^{k-1})^{-1}}
\big(\omega^{2/(2-\alpha)}\big)^{\sum_{k=1}^{i} (k-1)(p\omega^{k-1})^{-1}} \|f\|_{L^p(\Omega)}.
\end{equation*}
Since $\omega = n/(n-2) >1$, by letting $i \rightarrow \infty$, this implies
\begin{equation*}
\|f\|_{L^q (\Omega)}\lesssim\|f\|_{L^{\infty} (\Omega)}\lesssim C_\alpha^{\frac{n}{2p}} p^{\frac{n}{p(2-\alpha)}} \Big(\frac{n}{n-2} \Big)^{\frac{n(n-2)}{p(2-\alpha)}}\|f\|_{L^p (\Omega)}
\end{equation*}
as desired.

It remains to show \eqref{ineq1}.
Since $f \in H_0^1 (\Omega)$ is a positive part of the weak solution $u$, it follows that
\begin{equation} \label{fsol}
\int_{\Omega} \sum_{i,j=1}^{n} a_{ij}(x) f_{x_i} \phi_{x_j} dx+ \int_{\Omega}V(x) f(x) \phi (x) dx =  \int_{\Omega}\lambda f(x) \phi(x) dx
\end{equation}
for every real $\phi \in H_0^1(\Omega)$ supported on $\{x\in \mathbb{R}^n:u(x)>0\}$.
For $l>0$ and $m>0$, we set $\tilde f=f + l$, and let
\begin{equation*}
\tilde f_m = \begin{cases} l+m \quad \text{if} \quad f \geq m,\\ \tilde f \quad \text{if}\quad f<m.\end{cases}
\end{equation*}
We now consider the following test function
\begin{equation*}
\phi= \tilde f_m^{\beta} \tilde f - l^{\beta+1} \in H_{0}^1 (\Omega)
\end{equation*}
for $\beta \geq 0$.
We then compute
\begin{align*}
\phi_{x_j} &= \beta {\tilde f_m}^{\beta-1} (\tilde f_m)_{x_j}\tilde f + {\tilde f_m}^\beta \tilde f_{x_j}\\
&=\beta {\tilde f_m}^{\beta} (\tilde f_m)_{x_j} + {\tilde f_m}^\beta \tilde f_{x_j}
\end{align*}
using the fact that
\begin{equation}\label{concon}
 (\tilde f_m)_{x_j}=0\,\,\text{ in}\,\, \{x: f(x) \geq m\}\quad\text{and}\quad \tilde f_m = \tilde f\,\,\text{ in}\,\, \{x:f(x)<m\}.
\end{equation}

Substituting $\phi$ into \eqref{fsol} and using \eqref{concon} together with the trivial fact $f_{x_i} = \tilde f_{x_i}$,
the first term on the left-hand side of \eqref{fsol} is written as
\begin{align*}
\int_{\Omega} \sum_{i,j=1}^{n} & a_{ij}(x) f_{x_i} \phi_{x_j} dx \\
&= \beta \int_{\Omega} \tilde f_m^{\beta} \sum_{i,j=1}^n a_{ij}(x) (\tilde f_m)_{x_i} (\tilde f_m)_{x_j} dx + \int_{\Omega} {\tilde f_m}^\beta \sum_{i,j=1}^n a_{ij}(x) \tilde f_{x_i} \tilde f_{x_j} dx.
\end{align*}
Then it follows from the ellipticity \eqref{ellipticity} that
\begin{equation}\label{i}
\int_{\Omega} \sum_{i,j=1}^{n} a_{ij}(x) f_{x_i} \phi_{x_j} dx \geq \Lambda \beta \int_{\Omega} \tilde f_m^{\beta} |\nabla \tilde f_m|^2 dx + \Lambda \int_{\Omega} \tilde f_m^{\beta} |\nabla \tilde f|^2 dx.
\end{equation}
Combining \eqref{fsol} and \eqref{i}, we conclude that
\begin{align*}
\Lambda \beta \int_{\Omega} \tilde f_m^{\beta} |\nabla \tilde f_m|^2 dx &+ \Lambda \int_{\Omega} \tilde f_m^{\beta} |\nabla \tilde f|^2 dx \\
&\leq \int_{\Omega} -V f (\tilde f_m^{\beta} \tilde f - l^{\beta+1}) dx +  \int_{\Omega} \lambda f (\tilde f_m^{\beta} \tilde f -l^{\beta+1}) dx.
\end{align*}
Note here that
\begin{equation*}
|\nabla(\tilde f_m^{\beta/2}\tilde f)|^2  \leq 2 (1+\beta) \Big(\beta \tilde f_m^{\beta} |\nabla \tilde f_m|^2 + \tilde f_m^{\beta} |\nabla \tilde f|^2 \Big)
\end{equation*}
which follows by a direct computation together with \eqref{concon}.
We therefore get
\begin{align}\label{V}
\nonumber\int_{\Omega}  |\nabla &(\tilde f_m^{\beta/2} \tilde f)|^2 dx \\
\nonumber&\le \frac{2(1+\beta)}{\Lambda}\int_{\Omega} -Vf (\tilde f_m^{\beta}\tilde f-l^{\beta+1})dx + \frac{2\lambda(1+\beta)}{\Lambda} \int_{\Omega} f(\tilde f_m^{\beta}\tilde f-l^{\beta+1}) dx\\
&\le \frac{2(1+\beta)}{\Lambda}\int_{\Omega} |V| \tilde f_m^{\beta}\tilde f^2dx + \frac{2|\lambda|(1+\beta)}{\Lambda} \int_{\Omega} \tilde f_m^{\beta}\tilde f^2 dx.
\end{align}

To control the term involving the potential in \eqref{V}, we now use the so-called Fefferman-Phong inequality (\cite{F}),
\begin{equation}\label{FP}
\int_{\mathbb{R}^n}|g|^2 v(x)dx\leq C_n\|v\|_{\mathcal{L}^{2,r}} \int_{\mathbb{R}^n}|\nabla g|^2dx,
\end{equation}
where $C_n$ is a constant depending on $n$, and $1<r\leq n/2$.
(It is not valid for $r=1$ as remarked in \cite{D}.)
Applying this inequality along with H\"older's inequality, the first integral on the right-hand side of \eqref{V} is bounded as
\begin{align*}
\int_{\Omega} |V| \tilde f_m^{\beta} \tilde f^2 dx &\le \bigg( \int_{\Omega} |V|^{\frac{2}{\alpha}} \tilde f_m^{\beta} \tilde f^2 dx \bigg)^{\frac\alpha2} \bigg( \int_{\Omega} \tilde f_m^{\beta}
\tilde f^2 dx \bigg)^{\frac{2-\alpha}2} \\
&\le C_n\| |V|^{\frac{2}{\alpha}}\|_{\mathcal{L}^{2,\tilde{r}}}^{\frac{\alpha}{2}} \bigg( \int_{\Omega} |\nabla(\tilde f_m^{\beta/2} \tilde f)|^2 dx \bigg)^{\frac\alpha2} \bigg( \int_{\Omega}\tilde f_m^{\beta} \tilde f^2 dx \bigg)^{\frac{2-\alpha}2}
\end{align*}
for all $1<\tilde{r}\leq n/2$.
We note here that $\| |V|^{\frac{2}{\alpha}}\|_{\mathcal{L}^{2,\tilde{r}}}^{\frac{\alpha}{2}}=\|V\|_{\mathcal{L}^{\alpha,2\tilde{r}/\alpha}}$
and apply Young's inequality
\begin{equation} \label{Y}
ab \le \frac{\alpha}{2} (\varepsilon a)^{2/\alpha} + \frac{2-\alpha}{2} (\varepsilon^{-1} b)^{2/(2-\alpha)}
\end{equation}
with $0<\alpha <2$ and $\varepsilon >0$ to obtain
\begin{align} \label{e-ineq}
\nonumber
\int_{\Omega}& |V| \tilde f_m^{\beta} \tilde f^2 dx  \\
&\leq C_n\|V\|_{\mathcal{L}^{\alpha,2\tilde{r}/\alpha}} \bigg( \frac{\alpha}{2} \varepsilon^{2/\alpha} \int_{\Omega} |\nabla(\tilde f_m^{\beta/2}\tilde f)|^2 dx + \frac{2-\alpha}{2} \varepsilon^{{-2}/(2-\alpha)} \int_{\Omega}\tilde f_m^\beta \tilde f^2 dx \bigg).
\end{align}
By setting $r= 2\tilde{r}/\alpha$ (since $\widetilde{r}>1$, setting $r= 2\tilde{r}/\alpha$ determines the condition $r>2/\alpha$ in the theorem)
and taking $\varepsilon^{2/\alpha} =c(1+\beta)^{-1}$ with $c=\frac{\Lambda}{2\alpha C_n \|V\|_{\mathcal{L}^{\alpha,r}}}$
so that
$$C_n\|V\|_{\mathcal{L}^{\alpha,r}}\frac{\alpha}{2} \varepsilon^{2/\alpha}\frac{2(1+\beta)}{\Lambda}=1/2,$$
the gradient term in \eqref{e-ineq} can be absorbed into the left-hand side of \eqref{V}, as follows:
\begin{align}\label{bfGNS}
\nonumber\int_{\Omega} |\nabla (\tilde f_m^{\beta/2} \tilde f)|^2 dx &\leq
\,\,c^{-\frac2{2-\alpha}} \bigg( \frac{2-\alpha}{\alpha} \bigg) (\beta+1)^{\frac2{2-\alpha}}  \int_{\Omega} \tilde f_m^{\beta}\tilde f^2 dx \\
& + \frac{2|\lambda| (1+\beta)}{\Lambda} \int_{\Omega} \tilde f_m^{\beta} \tilde f^2 dx.
\end{align}

Finally, applying the Gagliardo-Nirenberg-Sobolev inequality (\cite{E}) to the left-hand side of \eqref{bfGNS}, we see
$$\bigg(\int_{\Omega}|\tilde f_m^{\beta/2}\tilde f|^{2\omega}dx\bigg)^{1/\omega}
\lesssim\int_{\Omega} |\nabla (\tilde f_m^{\beta/2} \tilde f)|^2 dx$$
with $\omega =n/(n-2)$.
Using the fact that $\tilde f_m \leq \tilde f$ and setting $\beta+2=\tau$, we therefore get
\begin{equation*}
\nonumber\bigg(\int_{\Omega}\tilde f_m^{\tau\omega} dx\bigg)^{1/\omega} \leq \,\,c^{-\frac2{2-\alpha}}\bigg( \frac{2-\alpha}{\alpha} \bigg) (\tau-1)^{\frac2{2-\alpha}}  \int_{\Omega} \tilde f^{\tau} dx + \frac{2|\lambda| (\tau-1)}{\Lambda} \int_{\Omega} \tilde f^{\tau} dx,
\end{equation*}
which implies the desired estimate
\begin{equation*}
\bigg( \int_{\Omega} {f}^{\tau\omega} dx \bigg)^{1/\omega} \lesssim \Big( 1+\alpha^{\frac{\alpha}{2-\alpha}}\Big(\frac{2C_n}{\Lambda}\|V\|_{\mathcal{L}^{\alpha, r}}\Big)^{\frac2{2-\alpha}} \Big) \tau^{\frac2{2-\alpha}} \int_{\Omega} f^{\tau} dx
\end{equation*}
by letting $m \rightarrow \infty$ and $l \rightarrow 0$.

\subsection{The case $p<2$}
From the case $p=2$, we have $\|f\|_{L^{\infty}(\Omega)} < \infty$ and
\begin{align}\label{p<2}
\|f\|_{L^\infty(\Omega)} &\lesssim C_\alpha^{\frac{n}{4}} 2^{\frac{n}{2(2-\alpha)}}\Big( \frac{n}{n-2}\Big)^{\frac{n(n-2)}{2(2-\alpha)}} \|f\|_{L^2(\Omega)} \nonumber \\
&\le C_\alpha^{\frac{n}{4}} 2^{\frac{n}{2(2-\alpha)}} \Big( \frac{n}{n-2}\Big)^{\frac{n(n-2)}{2(2-\alpha)}} \|f\|_{L^\infty(\Omega)}^{(2-p)/2} \|f\|_{L^p(\Omega)}^{p/2} \nonumber \\
&\le \frac{1}{2} \|f\|_{L^\infty(\Omega)} + C_\alpha^{\frac{n}{2p}}\frac{p}2\Big(\frac{1}{2-p}\Big)^{1-\frac2p} {2}^{\frac{n}{p (2-\alpha)}}\Big( \frac{n}{n-2}\Big)^{\frac{n(n-2)}{p(2-\alpha)}} \|f\|_{L^{p}(\Omega)}.
\end{align}
For the third inequality, we used here Young's inequality,
\begin{equation*}
ab \le \Big( 1- \frac{p}{2}\Big) (\epsilon a)^{\frac{2}{2-p}} + \frac{p}{2} (\epsilon^{-1} b)^{\frac{2}{p}}
\end{equation*}
with $\epsilon = ( \frac{1}{2-p})^{(2-p)/2}$,
\begin{equation*}
a=\|f\|_{L^{\infty}(\Omega)}^{(2-p)/2} \quad \textrm{and} \quad b=  C_\alpha^{\frac{n}{4}} 2^{\frac{n}{2(2-\alpha)}}
\Big( \frac{n}{n-2}\Big)^{\frac{n(n-2)}{2(2-\alpha)}} \|f\|_{L^p(\Omega)}^{p/2}.
\end{equation*}
By absorbing the first term on the right-hand side of \eqref{p<2} into the left-hand side, we conclude that
\begin{equation*}
\|f\|_{L^q (\Omega)}\lesssim\|f\|_{L^{\infty} (\Omega)}\lesssim C_\alpha^{\frac{n}{2p}} 2^{\frac{n}{p(2-\alpha)}} \Big(\frac{n}{n-2} \Big)^{\frac{n(n-2)}{p(2-\alpha)}}\|f\|_{L^p (\Omega)}
\end{equation*}
as desired.

\section{Concluding remarks} \label{sec3}
Finally we obtain a reverse H\"older inequality for $n=2$ by making use of a two-dimensional analogue (see \eqref{HS}) of the Fefferman-Phong inequality \eqref{FP}.

We first recall the function space $\mathcal{M}\log L$ introduced in \cite{A}, which is defined by
\begin{equation*}
	V \in  \mathcal{M} \log L \quad \Leftrightarrow \quad \sup_{0<t<|\Omega|}  \log \Big(\frac{|\Omega|}{t} \Big) \int_0^t V^\ast(s)ds < \infty .
\end{equation*}
For $t>0$, $V^\ast(t)= \inf \{s>0:d_V(s) \leq t\}$ is the decreasing rearrangement of $V$ with $d_V(s) = |\{x\in \Omega: |V(x)|> s\}|$ for $s>0$.   
Note that $L^p(\Omega) \subset L \log L \subset \mathcal{M} \log L$ (see \cite{A, BS}).

A two-dimensional version of Theorem \ref{thm1} is now stated in the following theorem.
Indeed, the constant $C_\alpha^{\frac{n}{2p}} \max\{p,2\}^{\frac{n}{p(2-\alpha)}} (\frac{n}{n-2})^{\frac{n(n-2)}{p(2-\alpha)}}$ in \eqref{RHI} boils down to the one in \eqref{RHI2} if $n\rightarrow2$.

\begin{thm}\label{thm2}
	Let $n=2$. Assume that $ u \in H_{0}^1 (\Omega) $ is a weak solution of the problem \eqref{maineq}
	with $\lambda\in\mathbb{R}$ and $V^{\sigma}\in \mathcal{M}\log L$ for some $\sigma>1$.
	Then we have
\begin{equation} \label{RHI2}
	\|u\|_{L^q (\Omega)}\leq C C_\sigma^{\frac{1}{p}} \max\{p,2\}^{\frac{\sigma}{p(\sigma-1)}} \|u\|_{L^p (\Omega)}
\end{equation}
	for all $q\geq p>0$. Here, $C$ is a constant depending on $\Lambda, \lambda, p, q$ and $\Omega$,
	and
	$$C_\sigma=1+\big(\frac{2}{\sigma}\big)^{\frac{1}{\sigma-1}}\Big(\frac{2\|V^{\sigma}\|_{\mathcal{M}\log L}}{\Lambda}\Big)^{\frac{\sigma}{\sigma-1}}.$$
\end{thm}

\begin{proof}
We briefly sketch the proof since it is an obvious modification of the one for Theorem \ref{thm1}. Indeed, it is enough to show the following inequality corresponding to \eqref{ineq1}: for all $\tau \ge 2$
\begin{equation}\label{ineq2}
	\|f\|_{L^{\tau \omega} (\Omega)} \lesssim C_\sigma^{1/\tau} \tau^{\frac{\sigma}{\tau(\sigma-1)}} \|f\|_{L^{\tau} (\Omega)}
\end{equation}
with $\omega= (2-\varepsilon)/\varepsilon$ for all arbitrarily small $\varepsilon>0$.
By iteration as before, we then get 
\begin{equation*}
	\|f\|_{L^q (\Omega)}\lesssim\|f\|_{L^{\infty} (\Omega)}\lesssim  C_\sigma^{\frac{2-\varepsilon}{2p(1-\varepsilon)}}p^{\frac{\sigma(2-\varepsilon)}{2p(\sigma-1)(1-\varepsilon)}} \Big(\frac{2-\varepsilon}{\varepsilon} \Big)^{\frac{\sigma\varepsilon(2-\varepsilon)}{2p(\sigma-1)(1-\varepsilon)^2}}\|f\|_{L^p (\Omega)}.
\end{equation*}
Letting $\varepsilon \rightarrow 0$ implies \eqref{RHI2} for $p\ge2$.
The case $p<2$ follows from the case $p=2$ in the same way as before.

It remains to show \eqref{ineq2}.
To control the term involving the potential in \eqref{V}, we use the following inequality (\cite{A})
\begin{equation}\label{HS}
	\int_{\Omega} |g|^2 v(x) dx \leq \frac{\|v\|_{\mathcal{M} \log L}}{\pi} \int_{\Omega} |\nabla g|^2 dx  
\end{equation}
instead of the Fefferman-Phong inequality \eqref{FP}.
Using \eqref{HS} along with H\"older's inequality, the first integral on the right side of \eqref{V} is indeed bounded as 
\begin{align*}
	\int_{\Omega} |V| \tilde f_m^{\beta} \tilde f^2 dx &\le \bigg( \int_{\Omega} |V|^{\sigma} \tilde f_m^{\beta} \tilde f^2 dx \bigg)^{\frac1\sigma} \bigg( \int_{\Omega} \tilde f_m^{\beta}
	\tilde f^2 dx \bigg)^{\frac{\sigma-1}\sigma} \\
	&\leq \| |V|^{\sigma}\|_{\mathcal{M} \log L} \bigg( \int_{\Omega} |\nabla(\tilde f_m^{\beta/2} \tilde f)|^2 dx \bigg)^{\frac1\sigma} \bigg( \int_{\Omega}\tilde f_m^{\beta} \tilde f^2 dx \bigg)^{\frac{\sigma-1}\sigma}
\end{align*}
for $\sigma>1$.
By using Young's inequality \eqref{Y} with $\alpha=2/\sigma$, we have
\begin{align} \label{e-ineq2}
	\nonumber
	\int_{\Omega}& |V| \tilde f_m^{\beta} \tilde f^2 dx  \\
	&\leq \||V|^\sigma\|_{\mathcal{M} \log L} \bigg( \frac{1}{\sigma} \varepsilon^{\sigma} \int_{\Omega} |\nabla(\tilde f_m^{\beta/2}\tilde f)|^2 dx + \frac{\sigma-1}{\sigma} \varepsilon^{{-\sigma}/(\sigma-1)} \int_{\Omega}\tilde f_m^\beta \tilde f^2 dx \bigg).
\end{align}
By taking $\varepsilon^{\sigma} =c(1+\beta)^{-1}$ with $c = \frac{ \sigma\Lambda}{4 \||V|^{\sigma}\|_{\mathcal{M}\log L}}$
so that
$$\||V|^\sigma\|_{\mathcal{M} \log L}\frac{1}{\sigma} \varepsilon^{\sigma}\frac{2(1+\beta)}{\Lambda}=\frac12,$$
the gradient term in \eqref{e-ineq2} can be absorbed into the left-hand side of \eqref{V}, as follows:
\begin{align}\label{bfGNS2}
	\nonumber\int_{\Omega} |\nabla (\tilde f_m^{\beta/2} \tilde f)|^2 dx &\leq
	\,\,c^{-\frac\sigma{\sigma-1}} (\sigma-1 ) (\beta+1)^{\frac\sigma{\sigma-1}}  \int_{\Omega} \tilde f_m^{\beta}\tilde f^2 dx \\
	& + \frac{2|\lambda| (1+\beta)}{\Lambda} \int_{\Omega} \tilde f_m^{\beta} \tilde f^2 dx.
\end{align}

Finally, we apply the Gagliardo-Nirenberg-Sobolev inequality (\cite{E}) after using H\"older's inequality in the left side of \eqref{bfGNS2} to obtain
\begin{equation}\label{ineq}
	\int_{\Omega} |\nabla (\tilde f_m^{\beta/2} \tilde f)|^2 dx \gtrsim \bigg(\int_{\Omega}   |\nabla(\tilde f_m^{\beta/2} \tilde f)|^{2-\varepsilon} dx \bigg)^{2/(2-\varepsilon)} \gtrsim	\bigg(\int_{\Omega}|\tilde f_m^{\beta/2}\tilde f|^{2\omega}dx\bigg)^{1/\omega}
\end{equation}
with $\omega= (2-\varepsilon)/\varepsilon$ for all arbitrarily small $\varepsilon>0$.
Using the fact that $\tilde f_m \leq \tilde f$ and setting $\beta+2 = \tau$, we therefore get 
\begin{equation*}
	\bigg( \int_{\Omega} {\tilde f_m}^{\tau\omega} dx \bigg)^{1/\omega} \lesssim \Big( 1+\big(\frac{2}{\sigma}\big)^{\frac{1}{\sigma-1}}\Big(\frac{2\||V|^{\sigma}\|_{\mathcal{M \log L}}}{\Lambda}\Big)^{\frac\sigma{\sigma-1}} \Big) \tau^{\frac\sigma{\sigma-1}} \int_{\Omega} \tilde f^{\tau} dx,
\end{equation*}
which implies the desired estimate \eqref{ineq2} by letting $m \rightarrow \infty$ and $l \rightarrow 0$.
\end{proof}


\end{document}